\documentclass[11pt,reqno]{amsart}

\usepackage[linktocpage=true]{hyperref}
\hypersetup{colorlinks,linkcolor=blue,urlcolor=blue,citecolor=blue}
\usepackage{amssymb}
\usepackage[cmtip,all]{xy}
\usepackage{enumerate} 
\usepackage{mathdots} 
\usepackage{blkarray} 
\usepackage{tikz}
\usetikzlibrary{cd,arrows}

\usepgflibrary{shapes.geometric}

\tikzstyle{V}=[draw, fill =black, circle, inner sep=0pt, minimum size=1.5pt]
\tikzstyle{C}=[draw, fill =white, circle, inner sep=0pt, minimum size=1.5pt]
\tikzstyle{over}=[draw=white,double=black,line width=2pt, double distance=.5pt]
\numberwithin{equation}{section}
\usepackage[margin = 0.9in]{geometry}
\theoremstyle{definition}

\newtheorem{theorem}{Theorem}[section]
\newtheorem{thm}[theorem]{Theorem}

\newtheorem{lemma}[theorem]{Lemma} 
\newtheorem{prop}[theorem]{Proposition}

\newtheorem{cor}[theorem]{Corollary}
\newtheorem{defn}[theorem]{Definition}
\newtheorem{definition}[theorem]{Definition}
\newtheorem{remark}[theorem]{Remark}

\newtheorem{ex}[theorem]{Example}

\newtheorem{rmk}[theorem]{Remark}

\def\<{\langle}
\def\>{\rangle}

\newcommand{\BD}{{I}^\textup{D}}

\newcommand{\cB}{\mathcal{B}}
\newcommand{\cBD}{\mathcal{B}^\textup{D}}
\newcommand{\CC}{\mathbb{C}}

\newcommand{\CP}{\mathbb{C}\mathrm{P}}

\newcommand{\D}{\textup{D}}

\newcommand{\tforall}{\textup{ for all }}

\newcommand{\GL}{\mathrm{GL}}

\newcommand{\inv}{^{-1}}
 
\newcommand{\II}{\textup{II}}
\newcommand{\III}{\textup{III}}

\newcommand{\ld}{\lambda}

\newcommand{\Sp}{\textup{Sp}}

\newcommand{\tif}{\textup{if }}

\newcommand{\ZZ}{\mathbb{Z}}

\allowdisplaybreaks 

\newcommand{\icupa}{
\begin{tikzpicture}[baseline={(0,-.3)}, scale = 0.8]
\draw (2.25,0) -- (0.75,0) -- (0.75,-.6) -- (2.25,-.6);
\draw (2.25,-.6) -- (2.25,0);\begin{footnotesize}
\node at (1,.2) {$1$};
\node at (1.75, -.3) {$a_2$};
\end{footnotesize}
\draw[thick] (1,0) -- (1, -.6);
\end{tikzpicture}
}

\newcommand{\icupb}{
\begin{tikzpicture}[baseline={(0,-.3)}, scale = 0.8]
\draw (2.25,0) -- (0.75,0) -- (0.75,-.6) -- (2.25,-.6);
\draw (2.25,-.6) -- (2.25,0);\begin{footnotesize}
\node at (1, -.3) {$\blacksquare$};
\node at (1,.2) {$1$};
\node at (1.75, -.3) {$a_2$};
\end{footnotesize}
\draw[thick] (1,0) -- (1, -.6);
\end{tikzpicture}
}

\newcommand{\icupc}{
\begin{tikzpicture}[baseline={(0,-.3)}, scale = 0.8]
\draw (3.25,0) -- (0.75,0) -- (0.75,-.9) -- (3.25,-.9);
\draw (3.25,-.9) -- (3.25,0);\begin{footnotesize}
\node at (2.5,.2) {$2t$};
\node at (1,.2) {$1$};
\node at (1.75, -.7) {$\blacksquare$};
\node at (1.75, -.3) {$a_1$};
\node at (2.8, -.3) {$a_2$};
\end{footnotesize}
\draw[thick] (1,0) .. controls +(0,-1) and +(0,-1) .. +(1.5,0);
\end{tikzpicture}
}

\newcommand{\icupd}{
\begin{tikzpicture}[baseline={(0,-.3)}, scale = 0.8]
\draw (3.25,0) -- (0.75,0) -- (0.75,-.9) -- (3.25,-.9);
\draw (3.25,-.9) -- (3.25,0);\begin{footnotesize}
\node at (2.5,.2) {$2t$};
\node at (1,.2) {$1$};
\node at (1.75, -.3) {$a_1$};
\node at (2.8, -.3) {$a_2$};
\end{footnotesize}
\draw[thick] (1,0) .. controls +(0,-1) and +(0,-1) .. +(1.5,0);
\end{tikzpicture}
}

\newcommand{\mcupea}{
\begin{tikzpicture}[baseline={(0,-.3)}, scale = 0.8]
\draw  (1.75,0) -- (-.25,0)  -- (-.25,-.9)--(1.75,-.9);
\draw (1.75,0) -- (1.75,-.9);
\begin{footnotesize}
\node at (0,.2) {$1$};
\node at (.5,.2) {$2$};
\node at (1,.2) {$3$};
\node at (1.5,.2) {$4$};
\end{footnotesize}
\draw[thick] (0,0) .. controls +(0,-.5) and +(0,-.5) .. +(.5,0);
\draw (1,0) -- (1,-.9);
\draw (1.5,0) -- (1.5,-.9);
\end{tikzpicture}
}

\newcommand{\mcupeb}{
\begin{tikzpicture}[baseline={(0,-.3)}, scale = 0.8]
\draw  (1.75,0) -- (-.25,0)  -- (-.25,-.9)--(1.75,-.9);
\draw (1.75,0) -- (1.75,-.9);
\begin{footnotesize}
\node at (0,.2) {$1$};
\node at (.5,.2) {$2$};
\node at (1,.2) {$3$};
\node at (1.5,.2) {$4$};
\end{footnotesize}
\draw[thick] (.5,0) .. controls +(0,-.5) and +(0,-.5) .. +(.5,0);
\draw (0,0) -- (0,-.9);
\draw (1.5,0) -- (1.5,-.9);
\end{tikzpicture}
}

\newcommand{\mcupec}{
\begin{tikzpicture}[baseline={(0,-.3)}, scale = 0.8]
\draw  (1.75,0) -- (-.25,0)  -- (-.25,-.9)--(1.75,-.9);
\draw (1.75,0) -- (1.75,-.9);
\begin{footnotesize}
\node at (0,.2) {$1$};
\node at (.5,.2) {$2$};
\node at (1,.2) {$3$};
\node at (1.5,.2) {$4$};
\end{footnotesize}
\draw[thick] (1,0) .. controls +(0,-.5) and +(0,-.5) .. +(.5,0);
\draw (0,0) -- (0,-.9);
\draw (0.5,0) -- (0.5,-.9);
\end{tikzpicture}
}

\newcommand{\mcuped}{
\begin{tikzpicture}[baseline={(0,-.3)}, scale = 0.8]
\draw  (1.75,0) -- (-.25,0)  -- (-.25,-.9)--(1.75,-.9);
\draw (1.75,0) -- (1.75,-.9);
\begin{footnotesize}
\node at (0,.2) {$1$};
\node at (.5,.2) {$2$};
\node at (1,.2) {$3$};
\node at (1.5,.2) {$4$};
\node at (0.5, -.5) {$\blacksquare$};
\node at (1.25, -.35) {$\blacksquare$};
\end{footnotesize}
\draw[thick] (1,0) .. controls +(0,-.5) and +(0,-.5) .. +(.5,0);
\draw (0,0) -- (0,-.9);
\draw (0.5,0) -- (0.5,-.9);
\end{tikzpicture}
}

\newcommand{\mcupee}{
\begin{tikzpicture}[baseline={(0,-.3)}, scale = 0.8]
\draw  (1.75,0) -- (-.25,0)  -- (-.25,-.9)--(1.75,-.9);
\draw (1.75,0) -- (1.75,-.9);
\begin{footnotesize}
\node at (0,.2) {$1$};
\node at (.5,.2) {$2$};
\node at (1,.2) {$3$};
\node at (1.5,.2) {$4$};
\node at (1.5, -.5) {$\blacksquare$};
\end{footnotesize}
\draw[thick] (0,0) .. controls +(0,-.5) and +(0,-.5) .. +(.5,0);
\draw (1,0) -- (1,-.9);
\draw (1.5,0) -- (1.5,-.9);
\end{tikzpicture}
}

\newcommand{\mcupef}{
\begin{tikzpicture}[baseline={(0,-.3)}, scale = 0.8]
\draw  (1.75,0) -- (-.25,0)  -- (-.25,-.9)--(1.75,-.9);
\draw (1.75,0) -- (1.75,-.9);
\begin{footnotesize}
\node at (0,.2) {$1$};
\node at (.5,.2) {$2$};
\node at (1,.2) {$3$};
\node at (1.5,.2) {$4$};
\node at (1.5, -.5) {$\blacksquare$};
\end{footnotesize}
\draw[thick] (0.5,0) .. controls +(0,-.5) and +(0,-.5) .. +(.5,0);
\draw (0,0) -- (0,-.9);
\draw (1.5,0) -- (1.5,-.9);
\end{tikzpicture}
}

\newcommand{\mcupeg}{
\begin{tikzpicture}[baseline={(0,-.3)}, scale = 0.8]
\draw  (1.75,0) -- (-.25,0)  -- (-.25,-.9)--(1.75,-.9);
\draw (1.75,0) -- (1.75,-.9);
\begin{footnotesize}
\node at (0,.2) {$1$};
\node at (.5,.2) {$2$};
\node at (1,.2) {$3$};
\node at (1.5,.2) {$4$};
\node at (0.5, -.5) {$\blacksquare$};
\end{footnotesize}
\draw[thick] (1,0) .. controls +(0,-.5) and +(0,-.5) .. +(.5,0);
\draw (0,0) -- (0,-.9);
\draw (0.5,0) -- (0.5,-.9);
\end{tikzpicture}
}

\newcommand{\mcupeh}{
\begin{tikzpicture}[baseline={(0,-.3)}, scale = 0.8]
\draw  (1.75,0) -- (-.25,0)  -- (-.25,-.9)--(1.75,-.9);
\draw (1.75,0) -- (1.75,-.9);
\begin{footnotesize}
\node at (0,.2) {$1$};
\node at (.5,.2) {$2$};
\node at (1,.2) {$3$};
\node at (1.5,.2) {$4$};
\node at (1.25, -.35) {$\blacksquare$};
\end{footnotesize}
\draw[thick] (1,0) .. controls +(0,-.5) and +(0,-.5) .. +(.5,0);
\draw (0,0) -- (0,-.9);
\draw (0.5,0) -- (0.5,-.9);
\end{tikzpicture}
}

\newenvironment{eq}{\begin{equation}}{\end{equation}}

\title{Irreducible components of two-row Springer fibers for all classical types}
\author[M.S. Im, C.-J. Lai, and A. Wilbert]{Mee Seong Im,  Chun-Ju Lai, and Arik Wilbert}
\address{Department of Mathematics, United States Naval Academy, Annapolis, MD 21402, USA}
    \email{meeseongim@gmail.com (Im)}
\address{Institute of Mathematics, Academia Sinica, Taipei 10617 Taiwan}
    \email{cjlai@gate.sinica.edu.tw (Lai)}
\address{Department of Mathematics and Statistics, University of South Alabama, Mobile, AL 36688, USA}
\email{wilbert@southalabama.edu (Wilbert)}
\begin{document}

\begin{abstract}
We give an explicit description of the irreducible components of two-row Springer fibers for all classical types using cup diagrams. Cup diagrams can be used to label the irreducible components of two-row Springer fibers. We use these diagrams to explicitly write down relations between the vector spaces of the flags contained in a given irreducible component. This generalizes results by Stroppel--Webster and Fung for type A to all classical types.  
\end{abstract}
 
\maketitle 


\section{Introduction}

Let $G$ be a connected complex reductive group with Lie algebra $\mathfrak{g}$. Given a nilpotent element $x\in\mathfrak{g}$, the corresponding Springer fiber $\mathcal B_x$ is defined as the variety of all Borel subgroups $B\subseteq G$ whose Lie algebra contains $x$. If $G=\GL_n$ (type A), the Springer fiber $\mathcal B_x$ can be realized more concretely as the variety of all complete flags $\{0\}=F_0\subseteq F_1\subseteq F_2\subseteq\ldots\subseteq F_n=\mathbb C^n$ satisfying $xF_i\subseteq F_{i-1}$ for all $1\leq i\leq n$. For the other classical groups $\mathrm{SO}_n$ (types B and D) and $\mathrm{Sp}_n$ (type C) one additionally imposes the condition $F_{n-i}^\perp=F_i$ for all $0\leq i\leq n$. Here, the orthogonal complement is taken with respect to some fixed nondegenerate symmetric (types B and D) or symplectic (type C) bilinear form compatible with the nilpotent $x\in\mathfrak{g}$. The varieties $\mathcal B_x$ naturally appear as the fibers under the Springer resolution and play an important role in classical geometric representation theory,~\cite{Spr76,Spr78}. Their geometric properties and combinatorial structure remain largely mysterious. In general, Springer fibers are not smooth and decompose into many equidimensional irreducible components.

In this article, we study Springer fibers for classical groups $G$, where the nilpotent element $x\in\mathfrak{g}$ has two Jordan blocks. Two-row Springer fibers are known to have some desirable geometric properties, e.g., the irreducible components are all smooth,~\cite{Fun03,FM10,ES16}. Moreover, they are known to have remarkable connections to representation theory and low-dimensional topology. 
For $G=\GL_n$, it was shown by Khovanov in ~\cite{Kho04} that the cohomology ring of a two-row Springer fiber is isomorphic to the center of the arc algebra appearing in his construction of link homology (see \cite{Kho02}). Arc algebras and their generalizations in~\cite{BS11,CK14} have shown intriguing connections to representation theory, e.g., they describe interesting categories of representations of the general Lie algebras,~\cite{BS11categoryO}, nonsemisimple walled Brauer algebras,~\cite{BS12brauer}, and the general linear Lie superalgebras,~\cite{BS12}. In~\cite{SW12}, Stroppel--Webster constructed the entire arc algebra (and their generalizations from~\cite{BS11,CK14}) as a convolution algebra using two-row Springer fibers.   
 
For two-row Springer fibers, the irreducible components can be classified in terms of so-called cup diagrams,~\cite{Fun03,ES16}. This diagrammatic tool does not only provide a direct connection to the combinatorics of the arc algebra, but it is particularly useful in writing down explicit relations characterizing all flags contained in a given irreducible component. For $\GL_n$, such an explicit description goes back to~\cite{SW12} based on~\cite{Fun03} (see Theorem~\ref{thm:mainA}). The goal of this article is to extend this result to two-row Springer fibers for all classical groups (see Theorem~\ref{thm:main1+2}). 

Studying the geometric, topological, and combinatorial properties of Springer fibers is an active field of research. Our work focusing on two-row Springer fibers is largely motivated by the fact that the center of the arc algebras of types B and D constructed in~\cite{ES15} are isomorphic to cohomology rings of two-row Springer fibers of types C and D,~\cite{ES16,SW16}. This relates these Springer fibers to representation theory of orthogonal Lie algebras,~\cite{ES15}, nonsemisimple Brauer algebras,~\cite{ES16_Brauer}, and orthosymplectic Lie superalgebras,~\cite{ES16super}.

Before this article was written, the only explicit algebro-geometric construction of irreducible components of two-row Springer fibers outside of type A required an intricate inductive procedure, see \cite[\S6]{ES16}, based on~\cite{Spa82, vL89}. Other than that, only topological models were available,~\cite{ES16,Wil18}.
 
This article is structured as follows. In Section 2 we provide the basic definitions, introduce the diagrammatic tools, state our main result, and discuss examples. It should be pointed out that our main theorem is not a straightforward generalization of the type A result from~\cite{SW12} and~\cite{Fun03}. In fact, due to the appearance of the bilinear form, quite a couple of new phenomena appear when generalizing the known results for type A to the other classical types. Our main theorem is proved in Section 3. The key idea is to prove that the explicit subvarieties we claim to be the irreducible components are iterated fiber bundles over $\CP^1$.    

\subsection*{Acknowledgments} We would like to thank Anthony Henderson and Catharina Stroppel for useful discussions and email exchanges. Some of the results contained in this article were originally available in the unpublished manuscript~\cite{ILW19}. The proof therein uses fixed-point subvarieties of certain Nakajima quiver varieties. We would like to thank Dongkwan Kim for pointing out that our main theorem can be obtained more efficiently without using the quiver variety. It has therefore been split off from the original preprint.        
 
C.-J. L. is partially supported by the MoST grant 109-2115-M-001-011-MY3, 2020 -- 2023.
\section{The main result}

In this section, we recall some basic facts about Springer fibers, state our main theorem, and provide examples.

\subsection{Nilpotent orbits and Springer fibers} \label{subsec:basic_facts_definitions} 
Let $n\in\mathbb Z_{>0}$ be a positive integer and let $V$ be an $n$-dimensional complex vector space. Let $G=\GL(V)$ be the group of all linear automorphisms of $V$. The Lie algebra $\mathfrak{g}=\mathfrak{gl}(V)$ of $G$ consists of all linear endomorphisms of $V$. 
Let $\mathcal N\subseteq\mathfrak{gl}(V)$ be the variety of all linear endomorphisms which are nilpotent. It is well known that the orbits under the conjugation-action of $G$ on $\mathcal N$ can be parametrized by partitions of $n$. Recall that a partition $\lambda$ of $n$, denoted by $\lambda\vdash n$, is a finite sequence $\lambda=(\lambda_1,\ldots,\lambda_r)\in\mathbb Z_{>0}^r$, $r\in\mathbb Z_{>0}$, of weakly decreasing positive integers $\lambda_1\geq\lambda_2\geq\ldots\geq\lambda_r$ summing up to $n$. Given a nilpotent orbit, the associated partition $\lambda$ encodes the Jordan type of the endomorphisms in that orbit, i.e., the parts of $\lambda$ encode the sizes of the different Jordan blocks.    

We fix $\epsilon\in\{\pm 1\}$ together with a non-degenerate bilinear form $\beta_\epsilon$ on $V$ satisfying $\beta_\epsilon(v,w)=\epsilon\beta_\epsilon(w,v)$ for all $v,w\in V$. Let $G_\epsilon$ be the group of all linear automorphisms of $V$ preserving $\beta_\epsilon$. Let $\mathfrak{g}_\epsilon$ be the Lie algebra of $G_\epsilon$, i.e., the subalgebra of $\mathfrak{g}$ consisting of all linear endomorphisms $x$ of $V$ satisfying $\beta_\epsilon(x(v),w)=-\beta_\epsilon(v,x(w))$ for all $v,w\in V$. We have isomorphisms $G_1\cong \textup{O}(n,\mathbb C)$ and $G_{-1}\cong \Sp(n,\mathbb C)$. Given $m\in\mathbb Z_{\geq 0}$, we refer to the groups $\GL(m,\mathbb C)$, $\textup{O}(2m+1,\mathbb C)$, $\Sp(2m,\mathbb C)$, and $\textup{O}(2m,\mathbb C)$ as the classical groups of types A, B, C, and D, respectively. 

Let $\mathcal P_\epsilon(n)$ be the subset of all partitions of $n$ defined by
\eq
\mathcal P_\epsilon(n)
=\{\lambda\vdash n \mid 
\#\{i \mid \lambda_i = j\} \in 2\ZZ
\tforall j \text{ with }(-1)^j=\epsilon
\}.
\endeq
The following result classifies the nilpotent orbits. 

\begin{prop}[{\cite{Ger61,Wi37}}]\label{prop:Gerstenhaber}
The orbits under the conjugation-action of $G_\epsilon$ on the variety $\mathcal N\cap\mathfrak{g}_\epsilon$ of all nilpotent elements in $\mathfrak{g}_\epsilon$ are in bijective correspondence with the partitions contained in $\mathcal P_\epsilon(n)$. The parts of the partition associated with the orbit of an endomorphism encode the sizes of the Jordan blocks in Jordan normal form.
\end{prop}


A complete flag in $V$ is a sequence $F_\bullet$ of nested vector subspaces $\{0\}=F_0\subsetneq F_1\subsetneq\ldots\subsetneq F_n=V$. We write $\mathcal B$ to denote the set of all flags in $V$. A complete flag $F_\bullet$ is called isotropic (with respect to $\beta_\epsilon$) if $F_{n-i}=F_i^\perp$ for all $i\in\{0,\ldots,n\}$. Here, the orthogonal complement is taken with respect to $\beta_\epsilon$. The set of all complete isotropic flags is denoted by $\mathcal B^{\beta_\epsilon}$. Both $\mathcal B$ and $\mathcal B^{\beta_\epsilon}$ can be identified with a homogeneous space, thereby giving these sets the structure of a smooth projective variety.


\begin{definition}
The Springer fiber $\mathcal B_x$ corresponding to $x\in\mathcal N$ is the subvariety of $\mathcal B$ consisting of all complete flags which satisfy the condition $xF_i\subseteq F_{i-1}$ for all $i\in\{1,\ldots,n\}$. The Springer fiber $\mathcal B^{\beta_\epsilon}_{x}$ associated with $\beta_\epsilon$ and $x\in\mathcal N\cap\mathfrak{g}_\epsilon$ is the subvariety of $\mathcal B^{\beta_\epsilon}$ consisting of all isotropic flags $F_\bullet$ which satisfy the conditions $xF_i \subseteq F_{i-1}$ for all $i\in\{1,\ldots,n\}$.
\end{definition}

\begin{remark}
We would like to point out that we work with the orthogonal group instead of the special orthogonal group in this article. This is mainly for technical reasons. In particular, the corresponding flag varieties and Springer fibers in type D decompose into two connected components (see e.g.~\cite[\S1.4]{vL89} for details). 
\end{remark}




In this article, we focus on two-row Springer fibers, i.e., the nilpotent endomorphism has Jordan type $\lambda=(n-k,k)$. Note that by the classification of nilpotent orbits in Proposition~\ref{prop:Gerstenhaber}, there is no two-row Springer fiber of type B. Moreover, we have the following result which reduces the problem of describing two-row Springer fibers for type C to the type D case.

\begin{prop}[{\cite[Theorem B]{Wil18}, \cite[Theorem 8.4.1]{Li19}}]\label{prop:iso_between_C_and_D}
Let $x\in\mathcal N\cap\mathfrak{g}_1$ be of Jordan type $(n-k,k)$ and $y\in\mathcal N\cap\mathfrak{g}_{-1}$ be of Jordan type $(n-k-1,k-1)$. Let $\overline{\mathcal B}^{\beta_1}_{x}$ be one of the two isomorphic connected components of the Springer fiber $\mathcal B^{\beta_1}_{x}$. Then we have an isomorphism of algebraic varieties $\overline{\mathcal B}^{\beta_1}_{x}\cong \mathcal B^{\beta_{-1}}_{y}$. 
\end{prop}

For the remainder of this article, we use the following explicit conventions for Springer fibers of types A and D: given a two-row partition $\ld=(n-k,k)$ of $n$, we define an $n$-dimensional complex vector space $V_\ld$ with (ordered) basis
\eq\label{eq:basisVld}
\{ e^\ld_1, e^\ld_2, \ldots, e^\ld_{n-k}, f^\ld_1, f_2^\ld, \ldots, f^\ld_k \}.
\endeq
Let $x_\lambda$ be the nilpotent endomorphism given by
\eq\label{eq:Jordanbasis}
e^\ld_{n-k} \mapsto e^\ld_{n-k-1} \mapsto \ldots \mapsto e^\ld_1 \mapsto 0,
\quad
f^\ld_{k} \mapsto \ldots \mapsto f^\ld_1 \mapsto 0.
\endeq
Usually, when there is no ambiguity, we omit the superscripts $\lambda$ in~\eqref{eq:basisVld} and~\eqref{eq:Jordanbasis} and we denote the endomorphism $x_{\ld}$ by $x$. The type A Springer fiber corresponding to these data is denoted by $\mathcal B_{n-k,k}$.
 
For the type D case, we assume that $n=2m$ is even, $\epsilon=1$, and $\ld \in \mathcal{P}_1$, i.e., $\ld$ is of the following form:
\eq\label{eq:partitionD}
\ld = (m,m),
\quad
\textup{or}
\quad
\ld = (n-k, k) \in (2\ZZ+1)^2.
\endeq
We also fix the non-degenerate symmetric bilinear form $\beta_\ld: V_\ld \times V_\ld \to \CC$, 
 whose associated matrix with respect to the ordered basis \eqref{eq:basisVld} is given by
\eq\label{eq:beta}
M^\ld=
\begin{cases}
\small
\begin{blockarray}{ *{8}{c} }
 & \{e^\ld_i\} & \{f^\ld_i\} \\
\begin{block}{ c @{\quad} ( @{\,} *{7}{c} @{\,} )}
\{e^\ld_i\}& 0 & J_m 
\\
\{f^\ld_i\}& J^t_m & 0
\\
\end{block}
\end{blockarray}~~~ 
\normalsize
&\:\:\tif n-k = k; 
\\
\small
\begin{blockarray}{ *{8}{c} }
 & \{e^\ld_i\} & \{f^\ld_i\} \\
\begin{block}{ c @{\quad} ( @{\,} *{7}{c} @{\,} )}
\{e^\ld_i\}& J_{n-k} & 0 
\\
\{f^\ld_i\}& 0 & J_k
\\
\end{block}
\end{blockarray}~~~
\normalsize
&\:\:\tif n-k > k,
\end{cases}
\quad
\textup{where}
\quad
J_i = 
\small
\begin{pmatrix}
&&&1
\\
&&-1&
\\
&\iddots&&
\\
(-1)^{i-1}
\end{pmatrix}
\normalsize.
\endeq
As above, we usually omit the superscripts in~\eqref{eq:beta} when the partition is fixed. The type D Springer fiber corresponding to the above data is denoted by $\cBD_{n-k,k}$.

\subsection{Irreducible components of type A two-row Springer fibers}

In~\cite{Spa76,Var79}, the irreducible components of type A Springer fibers were classified in terms of Young tableaux. For two-row Springer fibers, it is convenient to label the irreducible components by so-called cup diagrams,~\cite{Fun03}. For the rest of this subsection, we fix a two-row partition $\lambda=(n-k,k)$ of $n$.

\defn\label{def:cup_A}
We fix a rectangle in the plane together with $n$ evenly spaced vertices on the top horizontal edge of the rectangle. A cup diagram is a non-intersecting arrangement of cups and vertical rays inside the rectangle such that each vertex is connected to exactly one endpoint of either a cup or a ray. We write $I_{n-k,k}$ to denote the set of all cup diagrams on $n$ vertices with $k$ cups.
\enddefn

\begin{ex}\label{example:type_A_cup_diagrams}
The set $I_{3,1}$ consists of all cup diagrams on $n=4$ vertices with $k=1$ cup, i.e., 
\eq
\begin{array}{cc}
I_{3,1}=\left\{
a=\mcupea~,
\quad
b=\mcupeb~,
\quad
c=\mcupec
\right\}.
\end{array}
\endeq
\end{ex}

Given a cup diagram $a\in I_{n-k,k}$, we can define a subvariety of the Springer fiber $\mathcal B_{n-k,k}$ as follows. Here, we use the conventions introduced at the end of Subsection~\ref{subsec:basic_facts_definitions} omitting to indicate the partition $\lambda$.

\begin{defn}\label{def:main_type_A}
Given $a\in I_{n-k,k}$, define $K_{a}\subseteq \mathcal B_{n-k,k}$ as the subset consisting of all complete flags $F_\bullet$ satisfying the following conditions imposed by the diagram $a$:
\begin{enumerate}[(i)]
\item If vertices $i < j$ are connected by a cup, then 
\[
F_{j}=x^{-\frac{1}{2}(j-i+1)}F_{i-1}.
\]
\item If vertex $i$ is connected to a ray, then 
\[
F_i=\langle e_1,\ldots,e_{i-c(i)},f_1,\ldots,f_{c(i)}\rangle,
\]
where $c(i)$ is the number of cups to the left of $i$.
\end{enumerate}
\end{defn}

\begin{thm}[{\cite{SW12,Fun03}}]\label{thm:mainA}
The set $\{K_{a}~|~a\in I_{n-k,k}\}$ forms the complete list of irreducible components of the type A two-row Springer fiber $\mathcal B_{n-k,k}$. 
\end{thm}

\begin{ex}
For the cup diagrams from Example~\ref{example:type_A_cup_diagrams}, Definition~\ref{def:main_type_A} yields the following sets of flags
\begin{itemize}
\item $K_a=\big\{\left(\{0\}\subseteq F_1 \subseteq\<e_1,f_1\>\subseteq\<e_1,e_2,f_1\>\subseteq\<e_1,e_2,e_3,f_1\>=V\right)\big\}$,
\item $K_b=\big\{\left(\{0\}\subseteq\<e_1\>\subseteq F_2\subseteq\<e_1,e_2,f_1\>\subseteq\<e_1,e_2,e_3,f_1\>=V\right)\big\}$,
\item $K_c=\big\{\left(\{0\}\subseteq\<e_1\>\subseteq\<e_1,e_2\>\subseteq F_3\subseteq\<e_1,e_2,e_3,f_1\>=V\right)\big\}$.
\end{itemize}
By Theorem~\ref{thm:mainA}, these form a complete set of all irreducible component of $\mathcal B_{3,1}$. Note that each irreducible component is isomorphic to $\CP^1$. Using this explicit description of the irreducible components, it is straightforward to compute their intersections (see~\cite{Fun03} and~\cite{SW12} for details). As expected from a subregular nilpotent element, we get the following picture:
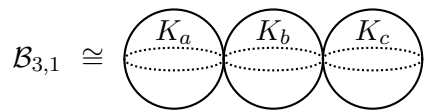
\[
\mathcal B_{3,1}\hspace{.4em}\cong\hspace{.4em}
\begin{tikzpicture}[baseline={(0,-.1)}, scale=.6]
\draw[thick] (0,0) circle (1.1);
\draw[thick] (2.2,0) circle (1.1);
\draw[thick] (4.4,0) circle (1.1);

\draw[densely dotted,thick] (0,0) ellipse (1.1 and 0.25);
\draw[densely dotted,thick] (2.2,0) ellipse (1.1 and 0.25);
\draw[densely dotted,thick] (4.4,0) ellipse (1.1 and 0.25);

\node at (2.2,.6) {$K_b$};
\node at (0,.6) {$K_a$};
\node at (4.4,.6) {$K_c$};
\end{tikzpicture}
\] 
\end{ex}

\subsection{Irreducible components of type D two-row Springer fibers}\label{sec:mcup}

For types B, C, and D, the irreducible components of Springer fibers can be classified in terms of so-called signed domino tableaux,~\cite{vL89,Spa82}. We now give an explicit characterization of the irreducible components of two-row Springer fibers for all classical types which generalizes the result for two-row Springer fibers of type A in Theorem~\ref{thm:mainA}. Recall that by Proposition~\ref{prop:iso_between_C_and_D} it suffices to treat the type D case. In this subsection, we assume that $n=2m$ is even and we fix a partition $\ld$ of the form as in~\eqref{eq:partitionD}.

In order to describe the irreducible components of type D two-row Springer fibers, we generalize Definition~\ref{def:cup_A} to also include so-called marked cup diagrams.

\defn\label{def:cup_D}
We call a cup or ray in a cup diagram accessible from the right if it can be connected to the right vertical edge of the rectangle by a path which does not intersect any other cup or ray. A marked cup diagram is a cup diagram in which any cup or ray that is accessible from the right is allowed to be decorated with a single marker (a small black square). We write $I_{n-k,k}^\D$ to denote the set of all marked cup diagrams on $m$ vertices with $\lfloor\frac{k}{2}\rfloor$ cups. This set decomposes as a disjoint union $I^\D_{n-k,k} = I_{n-k,k}^{\mathrm{even}} \sqcup  I_{n-k,k}^{\mathrm{odd}}$, where $ I_{n-k,k}^{\mathrm{even}}$ (resp.\  $I_{n-k,k}^{\mathrm{odd}}$) consists of all marked cup diagrams with an even (resp.\ odd) number of markers. 
\enddefn

\begin{ex}\label{ex:type_D_cup_diagrams}
The set $I_{5,3}^\D$ consists of the cup diagrams on $m=4$ vertices with $\lfloor\frac{3}{2}\rfloor=1$ cup, where 
\eq
\begin{array}{cc}
I_{5,3}^{\textup{even}}=\left\{
a=\mcupea~,
\quad
b=\mcupeb~,
\quad
c=\mcupec~,
\quad
d=\mcuped
\right\},
\\
I_{5,3}^{\mathrm{odd}} = \left\{
e=\mcupee~,
\quad
f=\mcupef~,
\quad
g=\mcupeg~,
\quad
h=\mcupeh
\right\}.
\end{array}
\endeq
\end{ex}

\rmk
Note that the marked cup diagrams only have $m=n/2$ vertices. This comes from the fact that marked cup diagrams in $I^\D_{n-k,k}$ are produced by folding certain symmetric cup diagrams on $n$ vertices in $I_{n-k,k}$, see~\cite{LS13}. The marked cup corresponds to a pair of unmarked cups crossing the axis of symmetry. More on this (including an algorithm on how to pass from ordinary cup diagrams to marked cup diagrams) can also be found in an unpublished manuscript~\cite{ILW19}. 
\endrmk

In type D, it is useful to distinguish between a two-row partition with equal and unequal parts. We start by considering the case in which $\ld = (m,m) = (k,k) \vdash n=2m$. We use the conventions for the type D Springer fiber $\cBD_{n-k,k}$ as discussed at the end of Subsection~\ref{subsec:basic_facts_definitions} omitting to indicate the partition $\lambda$ in the notation. 
\begin{defn}\label{def:main}
Given $a\in \BD_{k,k}$,
define $K_{a}\subseteq \cBD_{k,k}$ as the subset consisting of all complete flags $F_\bullet$ satisfying the following conditions imposed by the diagram $a$:
\begin{enumerate}[(i)]
\item If vertices $i < j$ are connected by a cup without a marker, then 
\[
F_{j}=x^{-\frac{1}{2}(j-i+1)}F_{i-1}.
\]
\item If vertices $i < j$ are connected by a marked cup, then 
\[
F_{i-1}+x^{\frac{1}{2}(j-i+1)}F_j = F_i \hspace{1em} \text{and} \hspace{1em} 
F_j^\perp = x^{-\frac{1}{2}(n-2j)}F_j .
\]
\item If vertex $i$ is connected to a marked ray, then 
\[
F_i=\langle e_1,\ldots,e_{\frac{1}{2}(i-1)},f_1,\ldots,f_{\frac{1}{2}(i+1)}\rangle.
\]
\item If vertex $i$ is connected to a ray without a marker, then 
\[
F_i=\langle e_1,\ldots,e_{\frac{1}{2}(i+1)},f_1,\ldots,f_{\frac{1}{2}(i-1)}\rangle.
\]
\end{enumerate}
\end{defn}

The key new phenomenon in type D when compared to type A is the appearance of the marked cup relation in Definition~\ref{def:main}(ii). We point out that this relation is in fact a pair of equations instead of a single equation. The purpose of the following example is to show in detail how to work with the equations for a marked cup. 

\begin{ex}\label{ex:UImU}
Let $F_\bullet \in K_{a}$, where $a \in \BD_{5,5}$ is the cup diagram below:
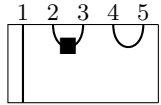
\[ 
\begin{tikzpicture}[baseline={(0,-.3)}, scale = 0.8]
\draw (3.75,0) -- (1.25,0)  -- (1.25,-1.3) --  (3.75,-1.3);
\draw (3.75,-1.3) -- (3.75,0);
\begin{footnotesize}
\node at (1.5,.2) {$1$};
\node at (2,.2) {$2$};
\node at (2.25, -.35) {$\blacksquare$};
\node at (2.5,.2) {$3$};
\node at (3,.2) {$4$};
\node at (3.5,.2) {$5$};
\end{footnotesize}
\draw[thick] (1.5,0) -- (1.5, -1.3);
\draw[thick] (2,0) .. controls +(0,-.5) and +(0,-.5) .. +(.5,0);
\draw[thick] (3,0) .. controls +(0,-.5) and +(0,-.5) .. +(.5,0);
\end{tikzpicture}
\]
Since vertex $1$ is connected to a ray without a marker, we have $F_1=\<f_1\>$ by relation (iv). Moreover, $F_2$ must be spanned by $f_1$ and some vector of the form $\lambda e_1+\mu f_2$, where $\lambda\neq 0$ or $\mu\neq 0$. In the following argument we assume that $\lambda$ and $\mu$ are both nonzero. (The cases where $\lambda=0$ or $\mu=0$ can be treated similarly.) Then $F_3$ is obtained by adding a vector of the form 
\begin{equation}\label{eq:marked_cup_example_eq1}
\alpha\left(\lambda e_1-\mu f_2\right)+\beta\left(\lambda e_2 + \mu f_3\right)
\end{equation} 
to the span of $f_1$ and $\lambda e_1+\mu f_2$. By the first relation in (ii), we have $F_1+xF_3=F_2$. Thus, $\beta\neq 0$. Otherwise, we would have $xF_3=0$, a contradiction to the required relation. Note that 
\begin{equation}\label{eq:marked_cup_example_eq2}
\alpha\left(\lambda e_3-\mu f_4\right)+\beta\left(\lambda e_4 + \mu f_5\right) \in x^{-2}F_3.
\end{equation}
Since $F_3^\perp=x^{-2}F_3$ by the second relation in (ii), the vectors in~\eqref{eq:marked_cup_example_eq1} and~\eqref{eq:marked_cup_example_eq2} are orthogonal and we get
\[
0=\beta\left(\alpha\left(\lambda e_1-\mu f_2\right)+\beta\left(\lambda e_2 + \mu f_3\right),\alpha\left(\lambda e_3-\mu f_4\right)+\beta_{5,5}\left(\lambda e_4 + \mu f_5\right)\right)=4\alpha\beta\lambda\mu
\] 
by the definition of the bilinear form $\beta_{5,5}$ in~\eqref{eq:beta}. Since $\beta,\lambda,\mu\neq 0$, we get $\alpha=0$. Hence, we have
\[
F_3=\<f_1,\lambda e_1+\mu f_2,\lambda e_2+\mu f_3\>,
\]
which is also true if $\lambda=0$ or $\mu=0$. By relation (i), we also get
\[
F_5=\<e_1,f_1,f_2,\lambda e_2+\mu f_3,\lambda e_3+\mu f_4\>.
\]
The space $F_4$ can be any four-dimensional space nested between $F_3$ and $F_5$. This describes all flags in $K_{a}$.
\end{ex}

Next, we deal with the second case in which  $n-k > k$.

\begin{defn}\label{def:main2}
Given $a\in \BD_{n-k,k}$,
define $K_{a}\subseteq \cBD_{n-k,k}$ as the subset consisting of all complete flags $F_\bullet$ satisfying (i)--(ii) of Definition~\ref{def:main} and the following conditions imposed by the diagram $a$:
\begin{enumerate}[(i)]
\setcounter{enumi}{2}
\item If vertex $i$ is connected to a marked ray, then 
\[
F_i=\langle e_1,\ldots,e_{i-c(i)-1},f_1,\ldots,f_{c(i)}, {f_{c(i)+1}+e_{i-c(i)}}\rangle.
\]
\item If vertex $i$ is connected to the rightmost ray without a marker, then 
\[
F_i=\langle e_1,\ldots,e_{i-c(i)-1},f_1,\ldots,f_{c(i)}, {f_{c(i)+1}-e_{i-c(i)}}\rangle.
\]
\item If vertex $i$ is connected to an unmarked ray that is not the rightmost, then 
\[
F_i=\langle e_1,\ldots,e_{i-c(i)},f_1,\ldots,f_{c(i)}\rangle.
\]
\end{enumerate}
\end{defn}

The following Theorem~\ref{thm:main1+2} is the main result and will be proved in Section~\ref{sec:pf}.

\begin{thm} \label{thm:main1+2}
The set $\{K_{a}~|~a\in \BD_{n-k,k}\}$ forms the complete list of irreducible components of the two-row Springer fiber $\cBD_{n-k,k}$ of type D. 
\end{thm}

\begin{ex}\label{ex:||U}
For the cup diagrams from Example~\ref{ex:type_D_cup_diagrams} without a marked cup, Definition~\ref{def:main2} yields the following sets of flags (we omit the orthogonal complements):
\begin{itemize}
\item $K_a=\big\{\left(\{0\}\subseteq F_1 \subseteq\<e_1,f_1\>\subseteq\<e_1,e_2,f_1\>\subseteq\<e_1,e_2,f_1,f_2-e_3\>\subseteq V\right)\big\}$,
\item $K_b=\big\{\left(\{0\}\subseteq\<e_1\>\subseteq F_2\subseteq\<e_1,e_2,f_1\>\subseteq\<e_1,e_2,f_1,f_2-e_3\>\subseteq V\right)\big\}$,
\item $K_c=\big\{\left(\{0\}\subseteq\<e_1\>\subseteq\<e_1,f_1-e_2\>\subseteq F_3\subseteq\<e_1,e_2,f_1,f_2-e_3\>\subseteq V\right)\big\}$,
\item $K_e=\big\{\left(\{0\}\subseteq F_1 \subseteq\<e_1,f_1\>\subseteq\<e_1,e_2,f_1\>\subseteq\<e_1,e_2,f_1,f_2+e_3\>\subseteq V\right)\big\}$,
\item $K_f=\big\{\left(\{0\}\subseteq\<e_1\>\subseteq F_2\subseteq\<e_1,e_2,f_1\>\subseteq\<e_1,e_2,f_1,f_2+e_3\>\subseteq V\right)\big\}$,
\item $K_g=\big\{\left(\{0\}\subseteq\<e_1\>\subseteq\<e_1,f_1+e_2\>\subseteq F_3\subseteq\<e_1,e_2,f_1,f_2+e_3\>\subseteq V\right)\big\}$.
\end{itemize}
In addition to that, the set $K_d$ consists of all flags where
\[
F_1=\<e_1\>\,\,\,,\,\,\,\,\,F_2=\<e_1,f_1+e_2\>\,\,\,,\,\,\,\,\,F_3=\<e_1,f_1+e_2,\lambda(f_1-e_2)+\mu(f_2+e_3)\>,
\]
\[
F_4=\<e_1,f_1+e_2,\lambda(f_1-e_2)+\mu(f_2+e_3),\lambda(f_2-e_3)+\mu(f_3+e_4)\>,
\]
and the set $K_h$ is given by all flags of the form
\[
F_1=\<e_1\>\,\,\,,\,\,\,\,\,F_2=\<e_1,f_1-e_2\>\,\,\,,\,\,\,\,\,F_3= \<e_1,f_1-e_2,\lambda(f_1+e_2)+\mu(f_2+e_3)\>,
\]
\[
F_4=\<e_1,f_1+e_2,\lambda(f_1-e_2)+\mu(f_2+e_3),\lambda(f_2+e_3)+\mu(f_3+e_4)\>,
\]
By Theorem~\ref{thm:main1+2}, the above sets form a complete set of all irreducible components of $\cBD_{5,3}$. Note that each irreducible component is isomorphic to $\CP^1$. Using the explicit description, it is straightforward to compute the intersections of the irreducible components. As expected from a subregular nilpotent element, we get the following picture:
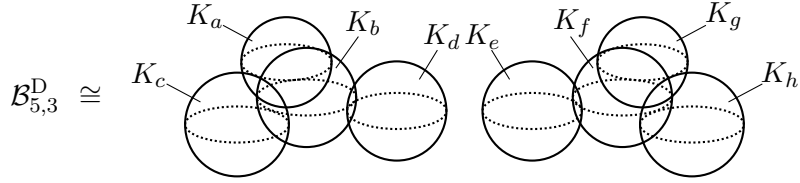
\[
\cBD_{5,3}\hspace{.4em}\cong\hspace{.4em}
\begin{tikzpicture}[baseline={(0,-.1)}, scale=.6]
\draw[thick] (-.445,.78) circle (1);
\draw[thick] (0,0) circle (1.1);
\draw[thick] (-1.54,-.6) circle (1.15);
\draw[thick] (2,-.3) circle (1.1);
\draw[densely dotted,thick] (0,0) ellipse (1.1 and 0.4);
\draw[densely dotted,thick] (2,-.3) ellipse (1.1 and 0.4);
\draw[densely dotted,thick] (-.445,.78) ellipse (1 and 0.4);
\draw[densely dotted,thick] (-1.54,-.6) ellipse (1.15 and 0.4);

\draw[thick] (7.445,.78) circle (1);
\draw[thick] (7,0) circle (1.1);
\draw[thick] (8.54,-.6) circle (1.15);
\draw[thick] (5,-.3) circle (1.1);
\draw[densely dotted,thick] (7,0) ellipse (1.1 and 0.4);
\draw[densely dotted,thick] (5,-.3) ellipse (1.1 and 0.4);
\draw[densely dotted,thick] (7.445,.78) ellipse (1 and 0.4);
\draw[densely dotted,thick] (8.54,-.6) ellipse (1.15 and 0.4);

\draw (.65,.65) -- (1.15,1.35);
\draw (2.5,.4) -- (2.9,1.05);
\draw (-2.35,-.1) -- (-3.1,.3);
\draw (-1.2,1.25) -- (-1.85,1.6);
\node at (1.25,1.65) {$K_b$};
\node at (3,1.35) {$K_d$};
\node at (-2.25,1.8) {$K_a$};
\node at (-3.5,.5) {$K_c$};

\draw (6.35,.65) -- (5.85,1.35);
\draw (4.5,.4) -- (4.1,1.05);
\draw (9.35,-.1) -- (10.1,.3);
\draw (8.2,1.25) -- (8.85,1.6);
\node at (5.9,1.65) {$K_f$};
\node at (3.9,1.35) {$K_e$};
\node at (9.25,1.8) {$K_g$};
\node at (10.5,.5) {$K_h$};
\end{tikzpicture}
\] 
Note that we have two connected components since we are working with the Springer fiber for the orthogonal group. Note that the parity of the markers on a cup diagram determines the connected component for the irreducible component labeled by the cup diagram.
\end{ex}

\section{Proof of the main theorem}\label{sec:pf}
In this section, we prove Theorem~\ref{thm:main1+2}. 
\subsection{Some small cases}
We first consider the case $n=2m \leq 4$.
\lemma\label{lem:n<=4}
Let $a\in \BD_{\ld}$ for some $\ld = (n-k,k)$, $n\leq 4$. Then $K_{a}$ is irreducible.
\endlemma
\proof
For each of the partitions $(1,1)$, $(2,2)$, and $(3,1)$, the set $a\in \BD_{\ld}$ contains two diagrams. It is easy to see that the subvariety $K_{a}$ is just a (geometric) point for each of the diagrams if the partition is either $(1,1)$ or $(3,1)$. For $(2,2)$, the two subvarieties are both isomorphic to $\CP^1$. In particular, they are irreducible.   
\endproof
\subsection{Iterated fiber bundles over $\CP^1$}

In order to prove the irreducibility of $K_{a}$ for $a \in B^\textup{D}_{n-k,k}$ in the general case, we recall the definition of an iterated fiber bundle over $\CP^1$.

\begin{definition}
A space $X$ is called an {\em iterated fiber bundle over $\CP^1$ of length $\ell$} if there exist spaces
$X=X_1, X_2, \ldots, X_\ell, X_{\ell+1} = \textup{pt}$ and maps $\pi_i:X_i \to \CP^1$ such that $\pi_i$ is a fiber bundle with typical fiber $X_{i+1}$ for $1\leq i \leq \ell$.
A point is considered an iterated bundle over  $\CP^1$ of length 0.
\end{definition}
The main goal of this section is to prove the following theorem.
\begin{thm}\label{thm:main3}
Let $a\in \BD_{n-k,k}$ be a marked cup diagram with $\ell$ cups. 
Then
$K_{a}$ is an iterated fiber bundle over $\CP^1$ of length $\ell$. 
\end{thm}

We have proved the special case $n\leq 4$ of Theorem~\ref{thm:main3} in Lemma~\ref{lem:n<=4}. 
Next, we will use an induction on $n=2m$. Since the inductive step is somewhat involved, we outline the idea here and provide details in the following subsections. 
Let $a\in \BD_{n-k,k}$ and let $K_{a}$ be as in Definitions~\ref{def:main} and~\ref{def:main2}. For $m \geq 3$, we have the following list of possible configurations of what the marked cup diagram $a \in \BD_{n-k,k}$ looks like locally at the first vertex:
\eq\label{eq:list}
\icupd~,
\quad
\icupc~,
\quad
\icupa~,
\quad
\icupb~
\quad
(1\leq t \leq \textstyle\lfloor\frac{m}{2}\rfloor).
\endeq
Here $a_2$ represents the (possibly empty) subdiagram of $a$ to the right of the cup or ray connected to vertex 1; while $a_1$ represents the (possibly empty) subdiagram of ${a}$ nested underneath the cup connected to vertex 1.
In order to prove that $K_{a}$ is an iterated fiber bundle over $\CP^1$, we distinguish three cases:

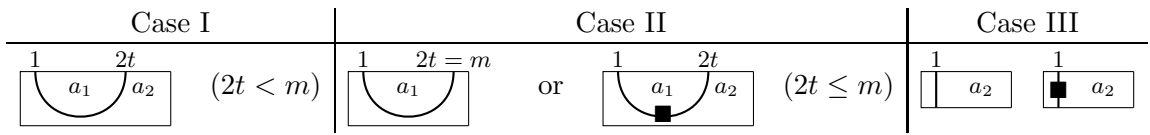
\[
\begin{tabular}{c|c|c}
Case I & Case II & Case III
\\
\hline
\icupd
\quad
$(2t<m)$
&
\begin{tikzpicture}[baseline={(0,-.3)}, scale = 0.8]
\draw (2.75,0) -- (0.75,0) -- (0.75,-.9) -- (2.75,-.9);
\draw (2.75,-.9) -- (2.75,0);\begin{footnotesize}
\node at (2.5,.2) {$2t=m$};
\node at (1,.2) {$1$};
\node at (1.75, -.3) {$a_1$};
\end{footnotesize}
\draw[thick] (1,0) .. controls +(0,-1) and +(0,-1) .. +(1.5,0);
\end{tikzpicture}
\quad
or
\quad
$\icupc 
\quad
(2t \leq m)$
&
$\icupa \quad \icupb$
\end{tabular}
\]
\begin{enumerate}
\item[Case I:] Vertex 1 is connected to $2t$ via an unmarked cup, where $2t < m$. In this case, $K_{a}$ is a trivial fiber bundle over $K_{b}$, i.e., $K_{a} \cong K_{b} \times K_{c}$, where  $K_{b}$ and $K_{c}$ are both iterated fiber bundles over $\CP^1$ of shorter lengths corresponding to two smaller cup diagrams $b$ and $c$ (see \eqref{eq:case1bc}). It then follows from~\cite[Lemma~8.11]{S12} that the product $K_{b} \times K_{c}$ is also an iterated fiber bundle over $\CP^1$. The details are provided in Subsection~\ref{sec:CaseI}.
\item[Case II:] Vertex 1 is connected to a cup, and is not Case I. In this case, $K_{a}$ is a nontrivial fiber bundle over $\CP^1$ with typical fiber $K_{c}$, where $c$ is a smaller cup diagram as in Definition~\eqref{def:ac}. By induction, $K_c$ is an iterated fiber bundle over $\CP^1$. The details are discussed in Subsection~\ref{sec:CaseII}. 
\item[Case III:] Vertex 1 is connected to a ray. In this case, $K_{a}$ is isomorphic to $K_{a_2}$, where $a_2$ is obtained from $a$ by removing the ray, and we can apply induction. More details can be found in Subsection~\ref{sec:CaseIII}.
\end{enumerate}

\subsection{Quadratic spaces}
In order to go through Cases I-III, we will need certain isomorphisms of {\em quadratic spaces}, i.e., vector spaces equipped with a symmetric bilinear form. 
For $\ld = (\ld_1, \ld_2) \vdash n$, we denote by $V_\ld$ the quadratic space $\CC^n \cong \<e_1^\ld, \ldots, e^\ld_{\ld_1}, f^\ld_1, \ldots, f^\ld_{\ld_2}\>$ equipped with the {bilinear form $\beta_\ld$ introduced in Subsection~\ref{subsec:basic_facts_definitions} (see \eqref{eq:basisVld}--\eqref{eq:beta}).}

For $\ld = (\ld_1,\ld_2), \mu = (\mu_1,\mu_2)$ such that $\ld_1 \geq \mu_1, \ld_2 \geq \mu_2$, we define a projection
\eq\label{def:P}
P^\ld_\mu: V_\ld \to V_\mu,
\quad
e^\ld_i \mapsto \begin{cases}
e^\mu_i &\tif 1\leq i \leq \mu_1;
\\
0 &\textup{otherwise},
\end{cases}
\quad
f^\ld_i \mapsto \begin{cases}
f^\mu_i &\tif 1\leq i \leq \mu_2;
\\
0 &\textup{otherwise}.
\end{cases}
\endeq
We also define the inverse map of its restriction on $\<e^\ld_i, f^\ld_j ~|~ 1\leq i \leq \mu_1, 1\leq j \leq \mu_2\>$ by
\eq
P^\mu_\ld: V_\mu \to V_\ld,
\quad
e^\mu_i \mapsto e^\ld_i,
\quad
f^\mu_j \mapsto f^\ld_j.
\endeq 
Now we consider an isotropic subspace $W \subseteq V_\ld$. 
Let  $\{x^\ld_i\}_{i\in I}$ and $\{x^\ld_j\}_{j\in J}$ be bases of $W$ and $W^\perp$, respectively, such that $I \subsetneq J$.
We write $\bar{x}^\ld_j = x^\ld_j + W$, and hence
$\{\bar{x}^\ld_j\}_{j \in J\setminus I}$ forms a basis of $W^\perp/W$. 
Denote the quotient map by $\psi = \psi(W)$ by
\eq\label{eq:psi}
\psi:W^\perp \to W^\perp/W.
\endeq
Moreover,  $W^\perp/W$ is a quadratic space equipped with the induced bilinear form $\bar{\beta}_\ld$ given by $\bar{\beta}_\ld(x+W,y+W) = \beta_\ld(x,y)$ for all $x,y \in W^\perp$.

In the following, we give explicit isomorphisms between certain quadratic spaces.
\begin{lemma}\label{lem:Q}
Let $\ld = (n-k,k)$.
\begin{enumerate}[(a)]
\item
Let $W$ be the subspace of $V_\ld$ spanned by $e_i^\ld,  f_i^\ld (1\leq i \leq t)$ such that $2t<m$.
Then 
\eq
W^\perp = \< e_i^\ld,  f_j^\ld ~|~ 1\leq i \leq n-k-t, 1\leq j \leq k-t \>.
\endeq
Moreover, there is a quadratic space isomorphism $Q^{\textup{I}}:W^\perp/W \to V_\nu$, for $\nu = \ld - (2t,2t)$, given by
\eq\label{eq:Qa}
\bar{e}^\ld_{t+i} \mapsto \begin{cases}
\sqrt{-1} e^\nu_i &\tif t \in 2\ZZ+1;
\\
e^\nu_i &\textup{otherwise},
\end{cases}
\quad
\bar{f}^\ld_{t+i} \mapsto \begin{cases}
\sqrt{-1} f^\nu_i &\tif t \in 2\ZZ+1;
\\
f^\nu_i &\textup{otherwise},
\end{cases}
\endeq
\item Assume that $\ld = (m,m)$ and $W = \<ce_1^\ld+ d f_1^\ld\>$ for some $(c,d) \in \CC^2\setminus \{(0,0)\}$. Then 
\[
W^\perp = \<c e_m^\ld+ d f_m^\ld, f^\ld_i ~|~ 1\leq i \leq m-1\>.
\]
Moreover, the assignments below all define quadratic space isomorphisms between $W^\perp/W$ and $V_\nu$, for $\nu = (m-1,m-1)$:
\begin{subequations}
\begin{align}
&Q^{\III}_1:&
\sqrt{-1} \bar{e}^\ld_{i+1} &\mapsto e^\nu_i,
&
\sqrt{-1} \bar{f}^\ld_{i} &\mapsto  f^\nu_i,
&\tif d=0;\label{eq:III-1}
\\
&Q^{\III}_2:&
\bar{f}^\ld_{i+1} &\mapsto f^\nu_i,
&
 \bar{e}^\ld_{i} &\mapsto  e^\nu_i,
&\tif c=0;\label{eq:III-2}
\\
&Q^{\II}_1:&
\sqrt{-1} (\overline{e^\ld_{i+1}+\textstyle\frac{d}{c} f^\ld_{i+1}}) &\mapsto e^\nu_i,
&
\sqrt{-1} \bar{f}^\ld_{i} &\mapsto  f^\nu_i,
&\tif c\neq 0, m\in2\ZZ;\label{eq:II-1}
\\
&Q^{\II}_2:&
(\overline{\textstyle\frac{c}{d} e^\ld_{i+1}+ f^\ld_{i+1}}) &\mapsto f^\nu_i,
&
\bar{e}^\ld_{i} &\mapsto  e^\nu_i,
&\tif d\neq 0, m\in2\ZZ.\label{eq:II-2}
\end{align}
\end{subequations}
\item Assume that $n-k > k$ and $W = \<e_1^\ld\>$. Then
\[
W^\perp = \<e_i^\ld, f^\ld_j ~|~ 1\leq i \leq n-k-1, 1\leq j \leq k\>.
\]
Moreover,  the assignments below all define quadratic space isomorphisms between $W^\perp/W$ and $V_\nu$, for $\nu = (n-k-2,k)$:
\begin{subequations}
\begin{align}
&Q^{\III}_3:&
\sqrt{\textstyle\frac{-1}{2}}(\overline{e^\ld_{i+1}+ f^\ld_{i}}) &\mapsto e^\nu_i,
&
\sqrt{\textstyle\frac{-1}{2}}(\overline{e^\ld_{i+1}- f^\ld_{i}}) &\mapsto  f^\nu_i,
&\tif n-k-2=k; \label{eq:III-3}
\\
&Q^{\III}_4:&
\sqrt{-1}\bar{e}^\ld_{i+1} &\mapsto e^\nu_i,
&
\sqrt{-1}\bar{f}^\ld_{i} &\mapsto  f^\nu_i,
&\tif  n-k-2>k. \label{eq:III-4}
\end{align}
\end{subequations}
\end{enumerate}
\end{lemma}
\proof
It follows from a direct computation using \eqref{eq:basisVld}--\eqref{eq:beta}.
For part (a), the form $\bar{\beta}_\ld$ is associated to the matrix obtained from $M_\ld$ by deleting columns and rows corresponding to $e^\ld_i, f^\ld_i$ for $1\leq i \leq t$.
Therefore, a naive projection $P^\ld_\nu$ does the job when $t$ is even; while in the odd case one needs to multiply the new basis elements by $\sqrt{-1}$ to make the forms to be compatible.

For part (b), we note first that $B_f = \{ ce^\ld_2+df^\ld_2, \ldots, ce^\ld_m+df^\ld_m, f^\ld_1, \ldots, f^\ld_{m-1}\}$ is an ordered basis of $W^\perp/W$ when $c \neq 0$; 
while $B_e = \{  e^\ld_1, \ldots, e^\ld_{m-1}, ce^\ld_2+df^\ld_2, \ldots, ce^\ld_m+df^\ld_m\}$ is an ordered basis of $W^\perp/W$ when $d \neq 0$.
In either case, the matrices associated to $\bar{\beta}_\ld$ with respect to $B_f, B_e$, respectively, are
\eq
\small
\begin{blockarray}{ *{8}{c} }
 & \{ce^\ld_{i+1}+df^\ld_{i+1}\} & \{f^\ld_i\} \\
\begin{block}{ c @{\quad} ( @{\,} *{7}{c} @{\,} )}
\{ce^\ld_{i+1}+df^\ld_{i+1}\}& A_{m-1} & -cJ_{m-1} 
\\
\{f^\ld_i\}& -cJ^t_{m-1} & 0
\\
\end{block}
\end{blockarray}
\quad
,
\quad
\begin{blockarray}{ *{8}{c} }
 & \{e^\ld_i\} & \{ce^\ld_{i+1}+df^\ld_{i+1}\} \\
\begin{block}{ c @{\quad} ( @{\,} *{7}{c} @{\,} )}
\{e^\ld_i\}& 0 & dJ_{m-1} 
\\
\{ce^\ld_{i+1}+df^\ld_{i+1}\}& dJ^t_{m-1} & A_{m-1}
\\
\end{block}
\end{blockarray}~~~ .
\normalsize
\endeq
Here $A_{m-1}$ is the zero matrix if $m$ is even; when $m=2r-1$ is odd, there only possible nonzero entry is
\eq
(A_{m-1})_{r,r} = \beta_\ld(ce^\ld_r+df^\ld_r, ce^\ld_r+df^\ld_r) = (-1)^{r-1}2cd.
\endeq
In other words, when $m$ is odd, $Q$ is an isomorphism of quadratic spaces if either $c=0$ or $d=0$.

For part (c), if $n-k-2 > k$, then the matrix of $\bar{\beta}_\ld$ is obtained from $M_\ld$ by deleting the rows and columns corresponding to $e_1^\ld$ and $e_m^\ld$, and hence one only needs to deal with the sign change of the upper left block. 
For the other case $n-k-2=k$, note first that $k$ must be odd so that \eqref{eq:partitionD} is satisfied. Next, we consider the ordered basis $B= \{e^\ld_2+f^\ld_1, \ldots, e^\ld_{k+1}+f^\ld_k, e^\ld_2-f^\ld_1, \ldots, e^\ld_{k+1}-f^\ld_k\}$. The matrix associated to $\bar{\beta}_\ld$ under $B$ is then
\eq
\small
\begin{blockarray}{ *{8}{c} }
 & \{e^\ld_{i+1}+f^\ld_i\} & \{e^\ld_{i+1}-f^\ld_i\} \\
\begin{block}{ c @{\quad} ( @{\,} *{7}{c} @{\,} )}
\{e^\ld_{i+1}+f^\ld_i\}& 0 & -2J_{k} 
\\
\{e^\ld_{i+1}-f^\ld_i\}& -2J^t_{k} &0
\\
\end{block}
\end{blockarray}~~~. 
\normalsize
\endeq 
We are done after a renormalization.
\endproof
For any $\mu = (\mu_1, \mu_2)$, let $x_\mu$ be the nilpotent operator given by
\eq
e^\mu_{\mu_1} \mapsto e^\mu_{\mu_1-1} \mapsto \ldots \mapsto e^\mu_1 \mapsto 0,
\quad
f^\mu_{\mu_2} \mapsto \ldots \mapsto f^\mu_1 \mapsto 0.
\endeq
Let W be one of the subspaces considered in cases (a)-(c) of Lemma~\ref{lem:Q}. 
Denote by
$\bar{x}_\ld$ be the induced nilpotent operator determined by
\eq\label{eq:barx}
\bar{e}^\ld_{i} \mapsto \bar{e}^\ld_{i-1},
\quad  
\bar{f}^\ld_{i} \mapsto \bar{f}^\ld_{i-1},
\endeq
For each such $W$, define an integer 
\eq
\ell = \begin{cases}
2t &\tif W = \<e_i, f_i ~|~ 1\leq i \leq t\>, 2t < m;
\\
1 &\textup{otherwise}.
\end{cases}
\endeq
We further denote by $\Omega = \Omega(W)$ the map
$\Omega: \cB^{\beta_\ld} \to \cB^{\beta_\nu}, F_\bullet \mapsto F''_\bullet$, where
\eq\label{eq:Omega}
F''_i = Q(F_{\ell+i}/W),
\quad (1\leq i \leq m-\ell)
\endeq
while $F''_{m} , \ldots, F''_{n-2\ell}$ are determined by the isotropy condition with respect to $\beta_\nu$.
\begin{cor}\label{cor:Q}
Using the notations in Lemma~\ref{lem:Q}, let $Q$ be one of the following quadratic space isomorphisms: $Q^{\textup{I}}, Q^{\II}_j (j=1,2), Q^{\III}_l (1\leq l \leq 4)$. Then
 $Q \bar{x}_\ld = x_\nu Q$.
Moreover, $\Omega(\cBD_{\bar{x}_\ld}) \subseteq \cBD_{x_\nu}$.
\end{cor}
\proof
The equality $Q \bar{x}_\ld = x_\nu Q$ is a direct consequence of Lemma~\ref{lem:Q} due to the explicit construction provided in \eqref{eq:Qa}--\eqref{eq:III-4}. 
For the latter statement, take any $F_\bullet \in \cB_{\bar{x}_\ld}$. Then we have
\[
 x_\nu Q(F_i/W)
=
Q \bar{x}_\ld (F_i/W)
\subseteq
Q (F_{i-1}/W).
\]
\endproof

In the next subsections, we use the following notation. For any cup diagram $a\in \BD_{n-k,k}$,
 denote the sets of all vertices connected to the left (resp.,\ right) endpoint of a marked cup in $a$ by $X_L^{a}$ (resp., $X_R^{a}$); while the sets of vertices connected to the left (resp., right) endpoint of an unmarked cup in $a$ are denoted by
$V_L^{a}$ (resp., $V_R^{a}$). If $i$ and $j$ are endpoints of a (marked) cup, define the endpoint-swapping bijection on $V_L^a \sqcup V_R^a \sqcup X_L^a \sqcup X_R^a$ by 
\eq\label{def:sig}
\sigma: 
i \mapsto j,
\quad
j\mapsto i.
\endeq
\subsection{Case I}\label{sec:CaseI}

In this section, we assume that
\eq\label{eq:caseI}
n =2m\geq 6,
\quad
\ld = (n-k,k),
\quad
a\in \BD_{\ld},
\quad
1 = \sigma(2t) \in V^{a}_L,
\quad
2t < m.
\endeq
In other words, we consider the marked cup diagram $a$ on $m \geq 3$ vertices with $\lfloor\frac{k}{2}\rfloor$ cups such that vertex 1 and $2t<m$ are connected by an unmarked cup, as below:
 \eq
a= \icupd
 \quad
 (2t<m). 
\endeq
We will show that $K_{a}$ is isomorphic to a trivial bundle $K_{{b}} \times K_{{c}}$ for some $b \in \BD_{\mu_1, \mu_2}$, $c \in \BD_{\nu_1, \nu_2}$ such that $\mu_2, \nu_2 < k$ so that the inductive hypothesis applies.
Note that if $2t=m$, the construction in this section will not work and has to be discussed in Section~\ref{sec:CaseII}
\begin{definition}\label{def:abc}
Assume that \eqref{eq:caseI} holds. We define two marked cup diagrams $b \in \BD_\mu, c \in \BD_\nu$ as follows:
\begin{enumerate}
\item $b$ is obtained by cropping $a$ on the first $2t$ vertices. Note that by construction, $b$ consists of only unmarked cups and hence $\mu =(2t,2t) \vdash 4t$.
\item $c$ is obtained by cropping $a$ on the last $m-2t$ vertices with a shift on the indices.
Note that $\nu =(n-2t-k,k-2t) \vdash n-4t$ by a simple bookkeeping on the number of vertices and cups.
\end{enumerate}
Namely, we have
\eq\label{eq:case1bc}
a= \icupd
\quad
\Rightarrow
b =
\begin{tikzpicture}[baseline={(0,-.3)}, scale = 0.8]
\draw (2.75,0) -- (0.75,0) -- (0.75,-.9) -- (2.75,-.9);
\draw (2.75,-.9) -- (2.75,0);\begin{footnotesize}
\node at (2.5,.2) {$2t$};
\node at (1,.2) {$1$};
\node at (1.75, -.3) {$a_1$};
\end{footnotesize}
\draw[thick] (1,0) .. controls +(0,-1) and +(0,-1) .. +(1.5,0);
\end{tikzpicture}
\quad
\textup{and}
\quad
c = 
\begin{tikzpicture}[baseline={(0,-.3)}, scale = 0.8]
\draw (2.25,0) -- (0.75,0) -- (0.75,-.9) -- (2.25,-.9);
\draw (2.25,-.9) -- (2.25,0);\begin{footnotesize}
\node at (2.25,.2) {$m-2t$};
\node at (1,.2) {$1$};
\node at (1.5, -.3) {$a_2$};
\end{footnotesize}
\end{tikzpicture}
\endeq
\end{definition}

\lemma\label{lem:Kc}
Let $a, c$ be defined as in Definition~\ref{def:abc}. 
Recall $\Omega$ from \eqref{eq:Omega} using $Q= Q^\textup{I}$ from \eqref{eq:Qa}.
Then $\Omega(K_{a}) \subseteq K_{c}$.
\endlemma
\proof
We first the the unmarked cup relations hold in $\Omega(K_{a})$.
In other words, for any $F_\bullet \in K_{a}$ and any pairs of vertices $(i,j)$ connected by an unmarked cup in $c$,
\eq\label{eq:newcup}
Q(F_{\ell+j}/W) = x_\nu^{\frac{-1}{2}(j-i+1)} Q(F_{\ell+i-1}/W),
\quad
\textup{for all}
\quad 
i \in V^{c}_L.
\endeq
Note that $(i-\ell, j-\ell)$ must be connected by an unmarked cup, and hence
\eq\label{eq:oldcup}
F_{j-\ell} = x_\ld^{\frac{-1}{2}(j-i+1)} F_{i-\ell},
\quad
\textup{for all}
\quad 
i \in V^{c}_L.
\endeq
Thus, \eqref{eq:newcup} follows from combining \eqref{eq:oldcup} and the former part of Corollary~\ref{cor:Q}. 

A ray connected to vertex $i$ in $c$ corresponds to the five flag conditions as in Definition~\ref{def:main}(iii)--(iv) and Definition~\ref{def:main2}(iii)--(v).
Since $b$ contains exactly $t$ unmarked cups, $F_{2t} = \<e_1, \ldots, e_t, f_1, \ldots, f_t\>$ and so 
$F''_i = Q(F_{\ell+i}/F_{2t})$. A case-by-case analysis shows that the new ray relations hold.
\endproof
\lemma\label{lem:Kb}
Let $a, b$ be defined as in Definition~\ref{def:abc}.
Then the map below is well-defined: 
\eq
\pi_{a,b}: K_{a} \to K_{b},
\quad
F_\bullet \mapsto F'_\bullet,
\endeq
where $F'_i = P^\ld_{\mu}(F_i)$ if $1 \leq i \leq 2t$, and that $F'_{2t+1} , \ldots, F'_{4t-1}$ are uniquely determined by $F'_1, \ldots, F'_{2t}$ under the isotropy condition with respect to $\beta_\mu$. 
\endlemma
\proof
We split the proof into three steps: firstly we show that $F'_\bullet$ is isotropic under $\beta_\mu$. Secondly, we show that $F'_\bullet$ sits inside the Springer fiber $\cBD_{x_\mu}$. 
Finally, we show that $F'_\bullet$ lies in the irreducible component $K_{b}$.

\begin{enumerate}[Step 1:]
\item
It suffices to show that $F'_{2t}$ is isotropic. Since the first $2t$ vertices in $a$ are all connected by unmarked cups, $F_{2t} = \<e^\ld_1, \ldots, e^\ld_t, f^\ld_1, \ldots, f^\ld_t\>$. Hence,
\eq
F'_{2t} = \<e^\mu_1, \ldots, e^\mu_t, f^\mu_1, \ldots, f^\mu_t\>.
\endeq
It then follows directly from \eqref{eq:beta} that $(F'_{2t})^\perp = F'_{2t}$.
\item
It suffices to check that $x_\mu F'_i \subseteq F'_{i-1}$ for all $i$. 
Note that a direct computation shows that
\eq\label{eq:Pxmu}
P_\mu^\ld x_\ld = x_\mu P_\mu^\ld.
\endeq
Hence, when $i\leq 2t$, we have
\eq
x_\mu F'_i  = x_\mu P^\ld_\mu(F_i) = P_\mu^\ld x_\ld (F_i) 
\subseteq 
P_\mu^\ld(F_{i-1}) = F'_{i-1}.
\endeq
Next, since now $x_\mu F'_{i+1} \subseteq F'_i$ for $i<2t$ , we have
$F'_{i+1} \subseteq x_\mu\inv F'_i$, and hence
\eq
F'_{4t-i-1}= (F'_{i+1})^\perp \supseteq (x_\mu\inv F'_i)^\perp = x_\mu F'_{4t-i}.
\endeq
\item  
We need to verify the conditions
\eq\label{eq:newcup}
P^\ld_\mu(F_{\sigma(i)}) = x_\mu^{\frac{-1}{2}(\sigma(i)-i+1)} (P^\ld_\mu(F_{i-1})),
\quad
\textup{for all}
\quad 
i \in V^{b}_L.
\endeq
Note that $V^{b}_L \subseteq V^{a}_L$, hence the unmarked cup relations in $a$ hold, i.e.,
\eq\label{eq:oldcup}
F_{\sigma(i)} = x_\ld^{\frac{-1}{2}(\sigma(i)-i+1)} F_{i-1},
\quad
\textup{for all}
\quad 
i \in V^{b}_L.
\endeq
Thus, \eqref{eq:newcup} follows from applying $P^\ld_\mu$ to \eqref{eq:oldcup},
thanks to \eqref{eq:Pxmu}. 
\end{enumerate}
\endproof
Finally, we are in a position to show that $K_{a}$ is irreducible.
\begin{prop}\label{prop:caseI}
Using the notations of Lemmas~\ref{lem:Kc}--\ref{lem:Kb}, the assignment
$F_\bullet \mapsto (\pi_{a,b}(F_\bullet), \Omega(F_\bullet))$ defines an isomorphism between $K_{a}$ and the trivial fiber bundle $K_{b} \times K_{c}$.
\end{prop}
\proof
It suffices to check that the given assignment is an isomorphism with inverse given by
$(F'_\bullet, F''_\bullet) \mapsto F_\bullet$, where
\eq
F_i = 
\begin{cases}
P^\mu_\ld(F'_i) &\tif 1\leq i \leq 2t;
\\
\psi\inv(Q\inv(F''_{i-2t})) &\tif 2t+1 \leq i \leq m;
\\
F_{n-i}^\perp &\tif m+1 \leq i \leq n.
\end{cases}
\endeq
That $F_\bullet \in K_{a}$ almost follows from construction, except for that we need to check if $F_{2t} \subset F_{2t+1}$, which follows from that $\psi\inv (Q\inv(F''_1))$ is a $2t+1$-dimensional space that contains $F_{2t}$.
\endproof
\subsection{Case II}\label{sec:CaseII}
In this section, we consider  one of the following diagrams:
\eq
\begin{tikzpicture}[baseline={(0,-.3)}, scale = 0.8]
\draw (2.75,0) -- (0.75,0) -- (0.75,-.9) -- (2.75,-.9);
\draw (2.75,-.9) -- (2.75,0);\begin{footnotesize}
\node at (2.5,.2) {$m=2t$};
\node at (1,.2) {$1$};
\node at (1.75, -.3) {$a_1$};
\end{footnotesize}
\draw[thick] (1,0) .. controls +(0,-1) and +(0,-1) .. +(1.5,0);
\end{tikzpicture}
\quad
\textup{or}
\quad
\icupc 
\quad
(2t \leq m)
\endeq
Note that for the latter subcase the marked cup $(1,2t)$ has to be accessible from the right vertical edge of the rectangle. So there are no rays in the subdiagram $a_2$. Moreover, there are no rays in the entire diagram $a$, and thus $\ld= (m,m)$.
That is, in this section we are working with the following assumptions:
\eq\label{eq:caseII}
n =2m\geq 6,
\quad
\ld = (m,m),
\quad
a\in \BD_{\ld},
\quad
\begin{array}{l}
\textup{either }
1 = \sigma(m) \in V^{a}_L
\\
\textup{or }
1 = \sigma(2t) \in X^{a}_L
\quad
(2t \leq m).
\end{array}
\endeq
We are going to show that $K_{a}$ is a fiber bundle over $\CP^1$ with typical fiber $K_{c}$ for some $c \in \BD_{\nu_1, \nu_2}$ such that $\nu_2 < k$ so that the inductive hypothesis applies.
\begin{definition}\label{def:ac}
Assume that \eqref{eq:caseII} holds. Define $c \in \BD_\nu$ by performing the following actions:
\begin{enumerate} 
\item Remove the cup connected to vertex 1 together with vertex 1. 
\item Connect a marked ray to vertex $\sigma(1)$.
\item Decrease all indices by one.
\end{enumerate}
Note that $\nu =(m-1,m-1)$ by a simple bookkeeping on the number of vertices and cups.
Namely, we have
\eq
a=
\begin{tikzpicture}[baseline={(0,-.3)}, scale = 0.8]
\draw (2.75,0) -- (0.75,0) -- (0.75,-.9) -- (2.75,-.9);
\draw (2.75,-.9) -- (2.75,0);\begin{footnotesize}
\node at (2.5,.2) {$m$};
\node at (1,.2) {$1$};
\node at (1.75, -.3) {$a_1$};
\end{footnotesize}
\draw[thick] (1,0) .. controls +(0,-1) and +(0,-1) .. +(1.5,0);
\end{tikzpicture}
\Rightarrow
c=
\begin{tikzpicture}[baseline={(0,-.3)}, scale = 0.8]
\draw (2.75,0) -- (1.25,0) -- (1.25,-.9) -- (2.75,-.9);
\draw (2.75,-.9) -- (2.75,0);\begin{footnotesize}
\node at (2.5,.2) {$m-1$};
\node at (2.5, -.65) {$\blacksquare$};
\node at (1.75, -.3) {$a_1$};
\end{footnotesize}
\draw[thick] (2.5,0) -- (2.5,-.9); 
\end{tikzpicture}
\quad
\textup{or}
\quad
a=
\icupc 
\Rightarrow
c=
\begin{tikzpicture}[baseline={(0,-.3)}, scale = 0.8]
\draw (3.25,0) -- (1.25,0) -- (1.25,-.9) -- (3.25,-.9);
\draw (3.25,-.9) -- (3.25,0);\begin{footnotesize}
\node at (2.5,.2) {$2t-1$};
\node at (2.5, -.65) {$\blacksquare$};
\node at (1.75, -.3) {$a_1$};
\node at (2.8, -.3) {$a_2$};
\end{footnotesize}
\draw[thick] (2.5,0) -- (2.5,-.9); 
\end{tikzpicture}
\endeq
\end{definition}
Define  $\pi:K_{a}\to \CP^1$ by
\eq
F_\bullet = ( 0 \subset F_1=\<\ld e_1+\mu f_1\> \subset \ldots \subset\CC^n) \mapsto [\ld:\mu].
\endeq
Now we check that $\pi$ is a fiber bundle with typical fiber $K_{c}$. For the local triviality condition we use the open covering $\CP^1 = U_1 \cup U_2$ where
\eq
U_1 = \{[1:\gamma] ~|~\gamma \in \CC\},
\quad
U_2 = \{[\gamma:1] ~|~\gamma \in \CC\}.
\endeq
\lemma\label{lem:KcII}
Let $a, c$ be defined as in Definition~\ref{def:ac}. 
Recall $\Omega$ from \eqref{eq:Omega} using $Q = Q^{\II}_j$, $j=1,2$ from Lemma~\ref{lem:Q}(b). 
Then the maps below are well-defined: 
\eq
\phi_j:\pi\inv(U_j) \to K_{c},
\quad
F_\bullet \mapsto \Omega(F_\bullet). 
\endeq
\endlemma
\proof
The cup relations in subdiagrams $a_1$ and $a_2$ can be verified similarly as in Lemma~\ref{lem:Kc}.
For the marked ray connected to vertex $2t-1$ in $c$, the corresponding relation, according to Definition~\ref{def:main}(iii) and Definition~\ref{def:main2}(iv), is
\eq\label{eq:newray}
F''_{2t-1} = 
\begin{cases}
\<e^\nu_1, \ldots, e^\nu_{t-1}, f^\nu_1, \ldots, f^\nu_t\>;
&\tif 2t=m;
\\
\<e^\nu_1, \ldots, e^\nu_{t-1}, f^\nu_1, \ldots, f^\nu_{t-1}, e^\nu_t+f^\nu_t\>
&\tif 2t<m.
\end{cases}
\endeq  
We may assume now $j=1$ since the other case can be proved by symmetry.
Take $F_\bullet \in \pi\inv(U_1)$ so that $F_1 = \<e_1^\ld + \gamma f^\ld_1\>$. The relation for the cup connecting $1$ and $2t$ is
\eq\label{eq:oldcupII}
\begin{cases}
F_{2t} =\<e^\ld_1, \ldots e^\ld_t, f^\ld_1, \ldots, f^\ld_t\>
&\tif 1 = \sigma(2t) \in V^{a}_L, 2t=m;
\\
x_\ld^t F_{2t} = F_1, 
&\tif 1 = \sigma(2t) \in X^{a}_L, 2t=m;
\\
x_\ld^t F_{2t} = F_1, 
\quad 
x^{m-2t} F_{2t}^\perp = F_{2t}
&\tif 1 = \sigma(2t) \in X^{a}_l, 2t < m.
\end{cases} 
\endeq
A direct application of Lemma~\ref{lem:Q}(b) shows that  \eqref{eq:oldcupII} leads to \eqref{eq:newray} in either case.
\endproof
\begin{prop}\label{prop:caseII}
Let $a, c$ be defined as in Definition~\ref{def:ac}. 
In both cases, the subvariety $K_{a}$ is a fiber bundle over $\CP^1$ with typical fiber $K_{c}$.
\end{prop}
\proof
For $j=1,2$, it suffices to check that the diagram below commutes:
\eq\label{eq:IIcomm}
\begin{tikzcd}
\pi\inv(U_j) \ar[r, "{(\pi,\phi_j)}","\simeq"'] \ar[d,"\pi_i"]
& U_j\times K_{c} \ar[ld, "\textup{proj}_1"]
\\
U_j
&
\end{tikzcd}
\endeq
It is routine to check that $(\pi,\phi_j)$ is an isomorphism with inverse given by
$([a:b], F''_\bullet) \mapsto F_\bullet$, where
\eq
F_i = 
\begin{cases}
\<ae_1^\ld+bf_1^\ld\> &\tif i=1;
\\
\psi\inv((Q^{\II}_j)\inv(F''_{i-1})) &\tif 2 \leq i \leq m;
\\
F_{n-i}^\perp &\tif m+1 \leq i \leq n.
\end{cases}
\endeq
\endproof
\subsection{Case III}\label{sec:CaseIII}
In this section, we assume that
\eq\label{eq:caseIII}
n =2m\geq 6,
\quad
\ld = (n-k,k),
\quad
a\in \BD_{\ld},
\quad
1 \not\in V^{a}_L \sqcup X^{a}_L
\endeq
That is, vertex 1 is connected to a ray.
We will discuss the following subcases:
\begin{enumerate}[~]
\item Case III--1: $a = \icupa$ in which $a_2$ contains no rays.
\item Case III--2: $a = \icupb$ in which $a_2$ contains no rays.
\item Case III--3: $a = \icupa$ in which $a_2$ contains exactly one ray.
\item Case III--4: $a = \icupa$ in which $a_2$ contains at least two rays.
\end{enumerate}
Let $c \in \BD_{\nu}$ be the marked cup diagram obtained from $a_2$ by decreasing all indices by one. 
We are going to show that $K_{a}$ is isomorphic to $K_{c}$ so that the inductive hypothesis applies since $\nu \vdash n-2 < n$.

For case III--$l$ ($1\leq l \leq 4$), by Lemma~\ref{lem:Q} we have a quadratic space isomorphism $Q^{\III}_l$. Below we summarize the data associated to each subcase:
\eq\label{eq:data}
\begin{array}{|c|c|c|c|c|}
\hline
\textup{Case} & \textup{III--}1& \textup{III--}2& \textup{III--}3& \textup{III--}4
\\
\hline
\ld & (m,m) &(m,m) & (m+1,m-1) & (n-k,k)
\\
W=F_1&\<e^\ld_1\> & \<f^\ld_1\> &\<e^\ld_1\>&\<e^\ld_1\>
\\
\nu & (m-1,m-1) & (m-1,m-1) & (m-1,m-1) & (n-k-2,k)
\\
Q^{\III}_l&\eqref{eq:III-1}&\eqref{eq:III-2}&\eqref{eq:III-3}&\eqref{eq:III-4}
\\
\hline
\end{array}
\endeq
\lemma\label{lem:KcIII}
Let $a, c$ be defined as in Definition~\ref{def:ac}. 
Recall $\Omega$ from \eqref{eq:Omega} using $Q = Q^{\III}_l$, $1\leq l \leq 4$ from Lemma~\ref{lem:Q}(b)--(c). 
Then the maps below are well-defined: 
\eq
\phi_l:  K_{a} \to K_{c},
\quad
F_\bullet \mapsto Q^{\III}_l(F_\bullet).
\endeq
\endlemma
\proof
The lemma can be proved using a routine case-by-case analysis in a similar way as Lemmas \ref{lem:Kc} and \ref{lem:KcII}.
\endproof
\begin{prop}\label{prop:caseIII}
Retain the notations of  Lemma~\ref{lem:KcIII}.
The map $\phi_l:  K_{a} \to K_{c}$ is an isomorphism.
\end{prop}
\proof
Recall the subspace $W$ from \eqref{eq:data} for each $l$.
It is routine to check that its inverse map is given by $F''_\bullet \mapsto F_\bullet$, where
\eq
F_i = 
\begin{cases}
W &\tif i=1;
\\
\psi\inv((Q^{\III}_l)\inv(F''_{i-1})) &\tif 2 \leq i \leq m;
\\
F_{n-i}^\perp &\tif m+1 \leq i \leq n.
\end{cases}
\endeq
\endproof

\subsection{Concluding the Proof}\label{subsection:conclude}
\proof[Proof of Theorem~\ref{thm:main3}]
By Lemma~\ref{lem:n<=4}, $K_{a}$ is an iterated fiber bundle over $\CP^1$ of length $\ell$ for $n \leq 4$.
For $n = 2m \geq 6$, an exhaustive list is given in \eqref{eq:list}, and can be divided into three cases.
In either case, $K_{a}$ is an iterated fiber bundle over $\CP^1$ of length $\ell$ thanks to Propositions~\ref{prop:caseI}, \ref{prop:caseII} and \ref{prop:caseIII}.
\endproof

\proof[Proof of Theorem~\ref{thm:main1+2}]
By Theorem~\ref{thm:main3} each $K_d$ for $d\in \BD_{n-k,k}$ is an irreducible component.
These components are distinct due to the constructive nature of Definitions~\ref{def:main} and \ref{def:main2}. 
Hence, since there exists a bijection between the irreducible components of the Springer fiber $\cBD_{n-k,k}$ and the set $\BD_{n-k,k}$ of marked cup diagrams (this follows by combining~\cite[Lemma~5.12]{ES16},~\cite[II.9.8]{Spa82} and~\cite[Lemmas~3.2.3, 3.3.3]{vL89}), we have described all the irreducible components.
\endproof
\bibliography{litlist-geom} \label{references}
\bibliographystyle{amsalpha}
 
\end{document}